\documentclass[10pt,letterpaper,reqno]{amsart}
\usepackage{amsmath,amssymb,amsthm}
\usepackage{hyperref}
\hypersetup{
  colorlinks=true,
  linkcolor=blue,
  citecolor=blue,
  urlcolor=blue
}

\newtheorem{theorem}{Theorem}
\newtheorem{lemma}[theorem]{Lemma}
\newtheorem{proposition}[theorem]{Proposition}
\newcommand{\C}{\mathbb C}
\newcommand{\R}{\mathbb R}
\newcommand{\T}{\mathbb T}
\newcommand{\Rea}{\operatorname{Re}}
\newcommand{\Id}{\operatorname{Id}} 
\newcommand{\norm}[1]{\left\|#1\right\|}

\title[Real alternative Daugavet property without the complex  Daugavet property]{Real alternative Daugavet property without Daugavet property in a  separable complex Banach space}
\author[Jinghao Huang]{Jinghao Huang}\address[Jinghao Huang]{Institute for Advanced Study in Mathematics, Harbin Institute of Technology, Harbin 150001, China \emph{E-mail~:} {\tt jinghao.huang@hit.edu.cn} }
\author[Ran Xu]{Ran Xu}\address[Ran Xu]{Institute for Advanced Study in Mathematics, Harbin Institute of Technology, Harbin 150001, China \emph{E-mail~:} {\tt xxurann@stu.hit.edu.cn} }

\thanks{ \footnotesize Jinghao Huang was supported by the NNSF of China (No. 12031004, 12301160, 12471134 and 12671159). }
\begin{document}

\begin{abstract}
We construct a separable complex Banach space whose underlying real space has
the alternative Daugavet property, although the complex space fails the
Daugavet property, thereby answering an open question posed by Kadets,
Mart\'in, and Pay\'a in 2006. 
\end{abstract}
\maketitle 
\section{Introduction}
A complex Banach space $X$ has the \emph{Daugavet property} (briefly, DP) if
\[
 \|\Id+T\|=1+\|T\|
\]
for every bounded complex-linear rank-one operator $T:X\to X$.
Here $\Id$ is the identity operator, and a rank-one operator has a
one-dimensional range. The underlying real space $X_{\R}$ consists of
the same vectors with the same norm, but only real scalar multiplication
is considered. The space $X_{\R}$  has the \emph{alternative Daugavet property} (briefly, ADP) if
\[
 \max\{\|\Id+R\|,\|\Id-R\|\}=1+\|R\|
\]
for every bounded real-linear rank-one operator $R:X_{\R}\to X_{\R}$; 
see \cite{adp}. Over a fixed scalar field, the DP implies
the ADP.  We refer to \cite{book} for background.

For real Banach spaces, the ADP does not imply the DP:
the real space $\ell_1$ has the ADP but does not have
the DP \cite{adp}.
We address the following question concerning complex
Banach spaces, posed by Kadets, Mart\'in, and Pay\'a
in their survey \cite[Problem~34]{survey}:
\begin{center}
\itshape
Must a complex Banach space $X$ have the DP
whenever $X_{\R}$ has the ADP?
\end{center}
The question is also listed as open in the monograph
draft by Kadets, Mart\'in, Rueda Zoca, and Werner
\cite[Question (12.5), p.~369]{book}.

The question has an affirmative answer for the complex
space $C(K,\C)$ of continuous functions on a nonempty
compact Hausdorff space $K$ and the complex space
$L_1(\mu;\C)$ of integrable functions for a nonzero
$\sigma$-finite measure $\mu$, equipped with their
usual norms, as follows by combining
the ADP characterizations of Mart\'in and Oikhberg
\cite[Theorem~3.4]{adp} with the classical DP
characterizations 
\cite[Section~2, Examples~(a) and~(b)]{Werner01}.

Our main result answers the question of Kadets,
Mart\'in, and Pay\'a in the negative, even for
separable complex Banach spaces.

For a real or complex Banach space $E$, let $B_E$ and $S_E$
denote its closed unit ball and unit sphere, respectively,
and let $E^*$ denote its dual over the corresponding scalar field. 
We use $\T=\{\omega\in\C:|\omega|=1\}$ and write $\Rea$ for the real part.

\begin{theorem}\label{main}
There exist a separable complex Banach space $X$, a vector $x\in S_X$,
and a functional $\phi\in S_{X^*}$ such that $\phi(x)=-1$,
$X_{\R}$ has the ADP, and the complex rank-one operator $T:X\to X$, defined by
$
 T(v)=\phi(v)x
$ for all $v\in X$,   satisfies
\[
 \|T\|=1,\qquad \|\Id+T\|<2=1+\|T\|.
\]
In particular, $X$ does not have the DP.
\end{theorem}
 
\section{A  finite-dimensional extension step for the counterexample}

We first establish an auxiliary inequality.

\begin{lemma}\label{signlemma}
For $a,q\in\C$ with $|a|,|q|\leq1$, there exists $\sigma\in\{-1,1\}$ such that
\begin{equation*}
 |a+\omega\sigma|+\Rea(\omega q)\leq\frac{3\sqrt3}{2}, 
 \qquad \forall \omega\in\T.
\end{equation*}
\end{lemma}
\begin{proof}
If $a=0$, choose $\sigma=1$. Then
\[
 |\omega\sigma|+\Rea(\omega q)\leq1+|q|\leq2, 
 \qquad \forall  \omega\in\T.
\] 
Assume that $a\neq 0$, and write $a=\eta |a|$, where  $\eta =a/|a|\in\T$. Choose $\sigma\in\{-1,1\}$ with $\Rea(\sigma\eta q)\leq0$. Put $b=\sigma\eta q$, and write $\omega=\sigma\eta e^{it}$, $-\pi\leq t\leq\pi$. Then,  
\[
|a+\omega\sigma|+\Rea(\omega q)=
 | |a| +e^{it}|+\Rea(be^{it}),\qquad |b|\leq1,\quad \Rea b\leq0.
\]

If $\cos t\geq0$, then
\[
|1+e^{it}|^2-| |a| +e^{it}|^2
=(1-|a|)(1+|a|+2\cos t)\ge0, 
\] 
which implies that $ | |a| +e^{it}|\leq|1+e^{it}|$. 
Since $ 1+e^{it}=e^{it/2}\bigl(e^{-it/2}+e^{it/2}\bigr) =
2e^{it/2}\cos(t/2), 
$ it follows that 
\[
 | |a| +e^{it}|\leq|1+e^{it}|=2\cos(t/2).
\] 
Since $\Rea b\le0$ and $|b|\le1$, 
it follows that for $|t|\le\pi/2$ we have
\[
\begin{aligned}
\Rea(be^{it})
&=(\Rea b)\cos t-(\operatorname{Im}b)\sin t
\le -(\operatorname{Im}b)\sin t\\
&\le |\operatorname{Im}b|\,|\sin t|
\le |b|\,|\sin t|
\le |\sin t|.
\end{aligned}
\]
The maximum of
$2\cos(t/2)+|\sin t|$ on $|t|\leq\pi/2$ is $3\sqrt3/2$,
attained at $|t|=\pi/3$.

If $\cos t\leq0$, then
$ 
 \bigl|\,|a|+e^{it}\bigr|^2
 =|a|^2+1+2|a|\cos t
 \leq|a|^2+1\leq2.
$ 
Moreover,
$ 
 \Rea(be^{it})\leq|be^{it}|=|b|\leq1.
$ 
Consequently,
\[
 \bigl|\,|a|+e^{it}\bigr|+\Rea(be^{it})
 \leq\sqrt2+1<\frac{3\sqrt3}{2}.
\] 
\end{proof}

The following proposition is the key step in the construction of our example.

\begin{proposition}\label{extension}
Let $E$ be a finite-dimensional complex Banach space,
and let $x\in S_E$ and $\phi\in S_{E^*}$ satisfy
\begin{equation}\label{inv}
 \phi(x)=-1,\qquad
 \|x+v\|_E+\Rea\phi(v)\leq\frac{13}{5}, 
 \quad \forall v\in B_E.
\end{equation} 
Then, for every $y\in B_E$ and $0<\delta<1$, there exist
a finite-dimensional complex Banach space $F$ containing
$E$ isometrically, an integer $N\geq2$, vectors
$u_1,\ldots,u_N\in S_F$, and a norm-one extension of
$\phi$ to $F$, still denoted by $\phi$, such that
\[
 \|x+v\|_F+\Rea\phi(v)\leq\frac{13}{5}, 
 \quad \forall v\in B_F,
\]
and
\begin{equation}\label{orthogonality}
\begin{aligned}
\frac1N\sum_{j=1}^N u_j
&=(1-\delta)y,\\[2pt]
\max\{\|z+u_j\|_F, & \|z-u_j\|_F\}
=2,
\qquad z\in S_E,\ 1\le j\le N.
\end{aligned}
\end{equation}
\end{proposition}
\begin{proof}
Choose an integer $N$ satisfying
\begin{equation}\label{eq:N-choice}
 N-1\geq\frac2\delta.
\end{equation} 
Let $F:=E\oplus\C^{N-1}$ be the algebraic direct sum.
We identify $E$ and $\C^{N-1}$ with
$E\oplus\{0\}$ and $\{0\}\oplus\C^{N-1}$,
respectively.
Every vector in $F$ then has a unique representation
$e+\xi$, where $e\in E$ and $\xi\in\C^{N-1}$.
Equip $F$ with the auxiliary norm
\[
 \|e+\xi\|_0:=\|e\|_E+\|\xi\|_2,
\qquad e\in E,\quad \xi\in\C^{N-1},
\]
where $\norm{\cdot}_2$ denotes the Euclidean norm. 
Take $v_1,\ldots,v_{N-1}$ to be the standard basis in $\C^{N-1}$, and put $v_N:=-\mathop{\sum}\limits_{j<N}v_j\in \C^{N-1}$. Define
\begin{align}\label{defuj}
 u_j=(1-\delta)y+v_j,
\qquad 1\le j\le N.
\end{align}
Since $\sum_{j=1}^N v_j =0$, 
it follows that  $ \frac1N\sum_{j=1}^N u_j
 =(1-\delta)y$.

Extend $\phi: E \to \mathbb{C}$ to $F$, still denoted by $\phi$, by setting 
\begin{align}\label{defphiextension}
     \phi: F=E\oplus\C^{N-1} &\notag \to \mathbb{C} \\ 
    e+\xi &  \mapsto  \phi(e)  \qquad(\forall e\in E,\ \xi\in\C^{N-1}).  
\end{align}
Set
\begin{equation}\label{eq:defq}
q:=\phi(u_j)\overset{\eqref{defuj}, \eqref{defphiextension}}{=}(1-\delta)\phi(y). 
\end{equation}
Since $
|\phi(y)|\le\|\phi\|_{E^*}\|y\|_E\le1$, 
it follows that 
$  |q|\leq1-\delta. $ 
For $1\leq j\leq N$ and $\omega\in\T$, set
\begin{equation}\label{eq:defw}
 w_{j,\omega}
 :=\frac{x+\omega u_j}
 {\frac{13}{5}-\Rea(\omega q)}.
\end{equation}
Let
\begin{equation}\label{generators}
 G:=B_E\cup\{u_j:1\leq j\leq N\}
 \cup\{w_{j,\omega}:1\leq j\leq N,  \omega\in\T\}.
\end{equation}
Define $B_F$  by 
\begin{equation}\label{eq:unitball-hull}
 B_F=\overline{\operatorname{conv}}^{\|\cdot\|_0}
 \{\rho g:g\in G,\ \rho\in\T\}.
\end{equation} 
Since $G$ is bounded with respect to $\norm{\cdot}_0$
and spans $F$, the set $B_F$ is closed and bounded
in $(F, \norm{\cdot}_0)$, as well as absolutely convex
and absorbing.
Consequently, the Minkowski functional
\[
\|v\|_F:=\inf\{\lambda>0:v\in\lambda B_F\},
\qquad v\in F,
\]
defines a complex norm on $F$ whose closed unit ball
is $B_F$; see \cite[Theorem~1.35 and the proof of
Theorem~1.39]{Rudin91}.


\medskip\noindent\textbf{Step 1: Prove 
$\|v\|_F=\|v\|_E$ for all $v\in E$ and $\|\phi\|_{F^*}=1$.}

Define \begin{align}\label{eq:projection}
 \Pi:F=E\oplus\C^{N-1}&\notag\to E\\
 e+\xi&\mapsto e   \qquad(\forall e\in E,\ \xi\in\C^{N-1}).
\end{align}
For $\omega\in\T$, we have
\begin{align}
 \|x+\omega(1-\delta)y\|_E+\Rea(\omega q)\notag & \overset{\eqref{eq:defq}}{=}
 \|(1-\delta)(x+\omega y)+\delta x\|_E
 +(1-\delta)\Rea(\omega\phi(y))\notag\\
 & \, \leq(1-\delta)
 \bigl(\|x+\omega y\|_E+\Rea(\phi( \omega y))\bigr)
 +\delta\|x\|_E\notag\\
 &\overset{\eqref{inv}}{\leq}
 \frac{13(1-\delta)}5+\delta
 =\frac{13-8\delta}{5}.
 \label{eq:old-estimate}
 \end{align}
For $e\in B_E$, $1\leq j\leq N$, and $\omega\in\T$, we have $\Pi (e) =e\in B_E$ and 
\begin{align*}
 \|\Pi (u_j)\|_E
& \overset{\eqref{defuj}}{=}(1-\delta)\|y\|_E
 \leq1-\delta,   \\ 
 \|\Pi ( w_{j,\omega})
 \|_E
 & \overset{\eqref{eq:defw}}{=}
 \frac{\|x+\omega(1-\delta)y\|_E}
 {13/5-\Rea(\omega q)}
 \overset{\eqref{eq:old-estimate}}{\leq}
 \frac{(13-8\delta)/5-\Rea(\omega q)}
 {13/5-\Rea(\omega q)}\leq1.
\end{align*}
Consequently,
\begin{equation}\label{eq:ball-projection}
 \Pi(G)\subseteq B_E,\qquad
 \Pi(B_F)\overset{\eqref{eq:unitball-hull}}{\subseteq} B_E.
\end{equation}
For every $v\in E$, it follows that 
\begin{equation}\label{eq:isometry}
 \|v\|_E\overset{ \eqref{eq:projection}}{=}
 \|\Pi (v)\|_E
 \overset{\eqref{eq:ball-projection}}{\leq}\|v\|_F
 \overset{
B_E\subseteq G\subseteq B_F
 }{\leq}\|v\|_E, \quad \text{i.e.,} \quad  \|v\|_E=\|v\|_F. 
\end{equation}

For every $v\in B_F$, write $v=e+\xi$, where
$e\in E$ and $\xi\in \C^{N-1}$.
Hence
\begin{align*}
|\phi(v)|
&=|\phi(e+\xi)|
\overset{\eqref{defphiextension}}{=}|\phi(e)|
\overset{\eqref{eq:projection}}{=}|\phi(\Pi(v))|\\
&\leq\|\phi\|_{E^*}\|\Pi(v)\|_E
=\|\Pi(v)\|_E
\overset{\eqref{eq:ball-projection}}{\leq}1,
\end{align*}
where we used $\|\phi\|_{E^*}=1$.
Therefore, 
$\|\phi\|_{F^*}\leq1$.
On the other hand, $\|x\|_F\overset{\eqref{eq:isometry}}{=}\|x\|_E=1$,
so $x\in B_F$ and
$ 
|\phi(x)|\overset{\eqref{inv}}{=}1.
$ 
Therefore,
\begin{equation*}
 \|\phi\|_{F^*}
=1.
\end{equation*}

\medskip\noindent\textbf{Step 2: Prove 
$\|u_j\|_F=1$ and
 $\max\{\|z+u_j\|_F,\|z-u_j\|_F\}
 =2
 $ 
for all $z\in S_E$,  $1\leq j\leq N$.}

Fix $z\in S_E$ and $j\in\{1,\ldots,N\}$.
By the Hahn--Banach theorem, choose $h\in S_{E^*}$ with
$h(z)=1$.  Since $|h(x)|\leq\|h\|_{E^*}\|x\|_E=1$
and $|q|\leq1-\delta<1$, Lemma~\ref{signlemma} gives
$\sigma\in\{-1,1\}$ such that
\begin{equation}\label{eq:chosen-sign}
 |h(x)+\omega\sigma|+\Rea(\omega q)
 \leq\frac{3\sqrt3}{2}, 
 \qquad \forall \omega\in\T.
\end{equation}

For $e\in E$ and $\lambda_1,\ldots,\lambda_N\in\C$,  set 
\begin{align}\label{eq:defhj}
 h_j:F&\notag\to\C\\
 e+\sum_{k=1}^N\lambda_k u_k
 &\mapsto h(e)+\lambda_j\sigma
 +\frac{N(1-\delta)h(y)-\sigma}{N-1}
 \sum_{\substack{1\leq k\leq N\\k\ne j}}\lambda_k.
\end{align} 
To verify that $h_j$ is well-defined, suppose that
\[
e+\sum_{k=1}^N\lambda_k u_k=0,
\qquad e\in E,\quad \lambda_1,\ldots,\lambda_N\in\C.
\]
 Since  $u_k=(1-\delta)y+v_k$, it follows that 
\begin{align}\label{linearhj}
\sum_{k=1}^N\lambda_kv_k=-
\left(e+(1-\delta)\left(\sum_{k=1}^N\lambda_k\right)y\right).
\end{align}
The vector
$e+(1-\delta)(\sum_{k=1}^N\lambda_k)y$ belongs to $E$,
whereas $\sum_{k=1}^N\lambda_kv_k$ belongs to
$\C^{N-1}$.
Since $(E\oplus \{0\} )  \cap (\{0\} \oplus \C^{N-1})=\{0\}$,
it follows that 
\[
\sum_{k=1}^N\lambda_kv_k=0.
\]
Using $v_N=-\sum_{k=1}^{N-1}v_k$, we obtain
\[
\sum_{k=1}^{N-1}(\lambda_k-\lambda_N)v_k=0.
\]
The linear independence of $v_1,\ldots,v_{N-1}$ gives
\begin{equation}\label{eq:relations}
\lambda_1=\cdots=\lambda_N.
\end{equation}
Consequently,
\[
e\overset{\eqref{linearhj}, \eqref{eq:relations}}{=}-N(1-\delta)\lambda_N y.
\]
It follows that
\begin{align*}
&h(e)+\lambda_j\sigma
+\frac{N(1-\delta)h(y)-\sigma}{N-1}
 \sum_{k\ne j}\lambda_k\\
&\quad
\overset{ \eqref{eq:relations}}{=}
-N(1-\delta)\lambda_Nh(y)+\lambda_N\sigma
+\lambda_N\bigl(N(1-\delta)h(y)-\sigma\bigr)
=0.
\end{align*}
Thus the formula defining $h_j$ is independent of
the representation of each vector in $F$, so $h_j$
is well-defined. 
Complex linearity follows directly
from the defining formula, and $h_j|_E=h$.

For $k\ne j$,
\begin{align}
 |h_j(u_k)-(1-\delta)h(y)|
 &\notag \overset{\eqref{eq:defhj}}{=}
 \frac{|(1-\delta)h(y)-\sigma|}{N-1}
 \leq\frac{(1-\delta)|h(y)|+1}{N-1}
 \\ & \,\, \leq\frac2{N-1}
 \overset{\eqref{eq:N-choice}}{\leq}\delta,
 \label{eq:b-close} 
 \end{align}
 which implies that 
 \begin{align}
 |h_j(u_k)|
  \overset{\eqref{eq:b-close}}{\leq} (1-\delta)|h(y)|+\delta
\le (1-\delta)+\delta=1.
 \label{eq:b-bound}
\end{align}
For $k\ne j$ and $\omega\in\T$,
\begin{align}
 &\qquad |h(x)+\omega h_j(u_k)|+\Rea(\omega q)\notag
\\ &\notag \,  =| 
h(x)+\omega(1-\delta)h(y)
 +\omega\bigl(h_j(u_k)-(1-\delta)h(y)\bigr)
 | +\Rea(\omega q) 
 \\
 & \,  \leq|h(x+\omega(1-\delta)y)|
       +|h_j(u_k)-(1-\delta)h(y)|+\Rea(\omega q)\notag\\
 & 
  \overset{ \eqref{eq:b-close}}
 {\leq} 
 \|x+\omega(1-\delta)y\|_E
       +\Rea(\omega q)+\delta \notag\\
 &
 \overset{\eqref{eq:old-estimate}}
 {\leq}
 \frac{13-8\delta}{5}+\delta
 \leq\frac{13}{5}.
 \label{eq:b-generator-bound}
\end{align}
For $k=j$ and $\omega\in\T$,
\begin{align}
 |h(x)+\omega h_j(u_j)|+\Rea(\omega q)
 \overset{\eqref{eq:defhj}}{=}
|h(x)+\omega\sigma|+\Rea(\omega q)
 \overset{\eqref{eq:chosen-sign}}{\leq}
 \frac{3\sqrt3}{2}<\frac{13}{5}.
 \label{eq:sign-generator-bound}
\end{align}

By \eqref{eq:defhj}, $h_j(u_j)=\sigma$.
Together with \eqref{eq:b-bound}, this gives
\[
 |h_j(u_k)|
 \overset{\eqref{eq:defhj}, \eqref{eq:b-bound}}{\leq}1, \qquad 1\le k\le N, 
\]
and, for $\omega\in\T$,
\begin{align*}
 |h_j(w_{k,\omega})|
 \overset{\eqref{eq:defw}}{=}
 \frac{|h(x)+\omega h_j(u_k)|}{13/5-\Rea(\omega q)}
 \overset{\eqref{eq:b-generator-bound},  
    \eqref{eq:sign-generator-bound}}{\leq}
    1.
\end{align*}
Together with $h_j|_E=h\in S_{E^*}$, these estimates give $\mathop{\sup}\limits_{g\in G}|h_j(g)|\leq1$.  Therefore, for every $v\in B_F$,  $ 
 |h_j(v)|\leq1$, i.e., 
\begin{equation}\label{eq:hj-norm}
 \|h_j\|_{F^*} \leq1.
\end{equation}
Thus,
\begin{align}\label{2z+sigmaujF}
2 & \notag \overset{h(z)=1}{=}|h(z) +\sigma^2|\overset{\eqref{eq:defhj}}{=}|h_j(z+\sigma u_j)|
 \overset{\eqref{eq:hj-norm}}{\leq}\|z+\sigma u_j\|_F
  \\ & \,\,\,\, \, \leq\|z\|_F+\|u_j\|_F
\overset{\eqref{eq:isometry}, \eqref{generators}, 
    \eqref{eq:unitball-hull}}{\leq} \|z\|_E+1 =2.
\end{align}
In particular,
\begin{equation}\label{eq:new-unit}
 \|u_j\|_F=1
\end{equation}
and
\[
2\overset{\eqref{2z+sigmaujF}}{=}
 \|z+\sigma u_j\|_F
 \leq\max_{\tau\in\{-1,1\}}\|z+\tau u_j\|_F
 \leq\|z\|_F+\|u_j\|_F
\overset{\eqref{eq:isometry},\eqref{eq:new-unit}}{=}2.
\]

\medskip\noindent\textbf{Step 3: Prove
$\|x+v\|_F+\Rea\phi(v)\leq\frac{13}{5}$
for all $v\in B_F$.}

For $e\in B_E$ and $\rho\in\T$,
\begin{align}\label{1.1}
 \|x+\rho e\|_F+\Rea\phi(\rho e)
 \overset{\eqref{eq:isometry}}{=}
 \|x+\rho e\|_E+\Rea\phi(\rho e)
 \overset{\eqref{inv}}{\leq}\frac{13}{5}.
\end{align} 
For $1\leq j\leq N$ and $\rho\in\T$,
\begin{align}\label{1.2}
 \|x+\rho u_j\|_F+\Rea\phi(\rho u_j)
 \overset{\eqref{eq:defw}, \eqref{eq:defq}}{=}
 \bigl(13/5-\Rea(\rho q)\bigr)\|w_{j,\rho}\|_F
 +\Rea(\rho q)
\overset{\eqref{eq:unitball-hull}}{\leq}\frac{13}{5}.
\end{align}
For $\rho,\omega\in\T$, we have
$|\omega q|=|q|\leq1-\delta$, 
\begin{align*}
 x+\rho w_{j,\omega}
 \overset{\eqref{eq:defw}}{=}
 \frac{\bigl(\frac{13}{5}-\Rea(\omega q)+\rho\bigr)x
       +\rho\omega u_j}
      {\frac{13}{5}-\Rea(\omega q)}
      \end{align*} 
  and \begin{align}\label{eq:rotated-generator}
 \phi(\rho w_{j,\omega})
\overset{\eqref{defphiextension}, \eqref{eq:defw}}{=}
\phi\left(
\frac{\rho(x+\omega u_j)}
{\frac{13}{5}-\Rea(\omega q)}
\right)
=\frac{\rho\bigl(\phi(x)+\omega\phi(u_j)\bigr)}
{\frac{13}{5}-\Rea(\omega q)}\overset{\eqref{inv},\eqref{eq:defq}}{=}
\frac{\rho(-1+\omega q)}
{\frac{13}{5}-\Rea(\omega q)}.
 \end{align}
 Since $|\omega q|\le1$, we have
\[
|\omega q-1|
\le\sqrt{2\bigl(1-\Rea(\omega q)\bigr)}.
\]
Consequently,
\begin{align*}
&\frac85\left(\frac{13}{5}-\Rea(\omega q)\right)
-2-|\omega q-1|\\
&\quad\ge
\frac{14}{25}+\frac85\bigl(1-\Rea(\omega q)\bigr)
-\sqrt{2\bigl(1-\Rea(\omega q)\bigr)}\\
&\quad=
\frac85\left(
\sqrt{1-\Rea(\omega q)}-\frac{5\sqrt2}{16}
\right)^2+\frac{99}{400}>0.
\end{align*}
It follows that
\begin{align}\label{Lemmaold}
1+\frac{2+|\omega q-1|}
          {13/5-\Rea(\omega q)}<1+\frac85=\frac{13}{5}.
\end{align}
Consequently,
\begin{align}\label{1.3}
 & \notag\qquad\,\,\, \,    \|x+\rho w_{j,\omega}\|_F
  +\Rea\phi(\rho w_{j,\omega})\\
 &\notag \quad  \,\overset{\eqref{eq:defw}, \eqref{eq:rotated-generator}}{=}
 \frac{\left\|
       \bigl(\frac{13}{5}-\Rea(\omega q)+\rho\bigr)x
       +\rho\omega u_j
       \right\|_F
       +\Rea\bigl(\rho(\omega q-1)\bigr)}
      {\frac{13}{5}-\Rea(\omega q)}\\
 &\notag\quad\,\,\,\,\, \,  \leq
 \frac{\bigl|\frac{13}{5}-\Rea(\omega q)+\rho\bigr|
       \|x\|_F+\|u_j\|_F+|\omega q-1|}
      {\frac{13}{5}-\Rea(\omega q)}\\
 & \quad \, \, \,\,\,   \overset{ \eqref{eq:new-unit}}{\leq}
 \frac{\frac{13}{5}-\Rea(\omega q)+2+|\omega q-1|}
      {\frac{13}{5}-\Rea(\omega q)}
      \overset{\eqref{Lemmaold}}{<}
 \frac{13}{5}.
\end{align}
Thus, combining \eqref{1.1}, \eqref{1.2} and \eqref{1.3}, we have
\begin{equation}\label{eq:generator-invariant}
 \|x+\rho g\|_F+\Rea\phi(\rho g)\le
 \frac{13}{5}, 
 \qquad \forall g\in G,\ \rho\in\T.
\end{equation}
The function $H:F\to\R$, defined by
\(
H(v):=\|x+v\|_F+\Rea\phi(v),
\) $v\in F,$
is continuous and convex.  By
$
 B_F
 \overset{\eqref{eq:unitball-hull}}{=}
 \overline{\operatorname{conv}}
 \{\rho g:g\in G,\ \rho\in\T\}
 $, 
 we have $$
 \|x+v\|_F+\Rea\phi(v)\leq\frac{13}{5}
$$ for all $v\in B_F$. 
\end{proof}

\section{Proof of Theorem \ref{main}}

We first recall the following geometric characterization of the ADP in terms of slices of the unit ball.  
\begin{proposition}\cite[Proposition 2.1]{adp}\label{Prop2.1}
A real  Banach space $(E, \norm{\cdot}_{E})$ has the alternative Daugavet property if and only if 
  for every $f\in S_{E^*}$, $z\in S_E$,
 and $\alpha,\varepsilon>0$, there exists $u\in B_E$ such that
 \[
  f(u)>1-\alpha
  \qquad\text{and}\qquad
  \max_{
  \tau\in\{-1,1\}
  }\|u+\tau  z\|_E\geq2-\varepsilon.
 \] 
\end{proposition}

\begin{proof}[Proof of Theorem~\ref{main}]

\textbf{Step 1: Construct a separable complex Banach space $X$,
$x\in S_X$, and $\phi\in S_{X^*}$ such that $\phi(x)=-1$ and
$\|x+v\|_X+\Rea\phi(v)\leq\frac{13}{5}$ for all $v\in B_X$.}

Let $E_1:=\C$ with its usual norm, and set $x:=1$.
Define $\phi:E_1\to\C$ by $\phi(a):=-a$.
Then
\[
 \|x\|_{E_1}=\|\phi\|_{E_1^*}=1,
 \qquad \phi(x)=-1.
\]
For every $v=s+it\in B_{E_1}$, we have
$s\geq-1$ and $|t|\leq1$.
Hence, 
\[
 \|x+v\|_{E_1}+\Rea\phi(v)
 =|1+s+it|-s
\overset{1+s\geq 0}{\leq }
1+s+|t|-s
 \leq2<\frac{13}{5}.
\]
Thus \eqref{inv} holds on $E_1$.

We now construct inductively finite-dimensional complex
Banach spaces
$ 
 E_1\subseteq E_2\subseteq E_3\subseteq\cdots,
$ 
with each inclusion being an isometry. 
Suppose that, for some integer $n\geq1$,  $E_n$ has been constructed and that
$x\in S_{E_n}$ and $\phi\in S_{E_n^*}$ satisfy
\eqref{inv}.
Since $B_{E_n}$ is compact, choose a finite set
$A_n\subseteq B_{E_n}$ such that
\begin{equation}\label{eq:net-cover}
 B_{E_n}\subseteq
 \bigcup_{y\in A_n}(y+2^{-n} B_{E_n}).
\end{equation}
Write $A_n=\{y_1,\ldots,y_m\}$ and set $F_0:=E_n$.
We construct finite-dimensional complex Banach spaces
$F_1,\ldots,F_m$ as follows.
For each $i=1,\ldots,m$, apply Proposition~\ref{extension}
with $E=F_{i-1}$, $y=y_i$, and $\delta=2^{-n}$.
This gives a space $F_i$ containing $F_{i-1}$
isometrically and a norm-one extension of $\phi$,
still denoted by $\phi$, satisfying
\[
 \phi(x)=-1,\qquad
 \|x+v\|_{F_i}+\Rea\phi(v)\leq\frac{13}{5}, 
 \quad \forall v\in B_{F_i}.
\]
Finally, set $E_{n+1}:=F_m$. Thus
\[
 E_n=F_0\subseteq F_1\subseteq\cdots
 \subseteq F_m=E_{n+1},
\]
with all inclusions isometric. In particular, $E_{n+1}$ is finite-dimensional,
\begin{equation}\label{eq:stage-embeddings}
 \|v\|_{E_{n+1}}=\|v\|_{E_n}, \quad \forall v\in E_n,
 \qquad \|\phi\|_{E_{n+1}^*}=1,
\end{equation}
and
\begin{equation}\label{eq:stage-invariant}
 \phi(x)=-1,\qquad
 \|x+v\|_{E_{n+1}}+\Rea\phi(v)\leq\frac{13}{5}, 
 \quad \forall v\in B_{E_{n+1}}.
\end{equation}
For each $y\in A_n$, the construction in Proposition~\ref{extension}  yields an integer
$N\geq2$ and unit vectors $u_1,\ldots,u_N\in S_{E_{n+1}}$
such that
\begin{align*}
 \frac1N\sum_{j=1}^N u_j
 &\overset{\eqref{orthogonality}}{=}
 (1-2^{-n})y,
 \quad 
 \max_{\tau\in\{-1,1\}}\|z+\tau u_j\|_{E_{n+1}}
 \overset{\eqref{orthogonality}}{=} 2, \quad \forall z\in S_{E_n},\ 1\leq j\leq N.
\end{align*}

The isometric inclusions define a norm on
$\bigcup_{n\geq1}E_n$ by
\[
 \|v\|:=\|v\|_{E_n},\qquad \forall v\in E_n.
\]
This is well defined by \eqref{eq:stage-embeddings}.
Let $X$ be the completion of $\bigcup_{n\geq1}E_n$
with respect to this norm.
Then
\begin{equation}\label{eq:stage-isometry}
 \|v\|_X=\|v\|_{E_n}, 
 \qquad \forall n\geq1,\ v\in E_n.
\end{equation}
Since each $E_n$ is finite-dimensional and hence separable, choose a countable
dense subset $D_n\subseteq E_n$.
The set $\bigcup_{n\geq1}D_n$ is countable and 
\[
 \overline{\bigcup_{n\geq1}D_n}^{\,X}
 =\overline{\bigcup_{n\geq1}E_n}^{\,X}=X.
\]
Hence $X$ is separable.
Moreover, \begin{equation}\label{eq:ball-density}
 \overline{\bigcup_{n\geq 1} B_{E_n}}^{\,X}=B_X. 
\end{equation}

The complex-linear map
$ 
 \phi:\bigcup_{n\geq1}E_n\rightarrow\C
$  is well defined and 
 is bounded
with norm one. It therefore extends uniquely to a continuous
complex-linear functional on $X$, still denoted by $\phi$.
The functional $\phi: X \to \mathbb{C}$  satisfies
\begin{equation}\label{eq:limit-phi}
 \phi(x)=-1,\qquad
 \|\phi\|_{X^*}=1.
\end{equation}
By \eqref{eq:stage-invariant}, \eqref{eq:stage-isometry},
and the density assertion \eqref{eq:ball-density},
continuity yields 
\begin{equation}\label{eq:limit-invariant}
 \|x+u\|_X+\Rea\phi(u)
 \leq\frac{13}{5}, \qquad \forall u\in B_X.
\end{equation}
\medskip\noindent\textbf{Step 2: $X_{\R}$ has the ADP.}

Let $z\in S_X$, let $f\in S_{(X_{\R})^*}$ be a real-valued
real-linear functional, and let $\alpha,\varepsilon>0$.
Set
\begin{equation}\label{eq:beta}
 \beta:=\frac{\min\{\alpha,1\}}4,
 \quad   \text{then} \quad 0<\beta\leq\frac14,\qquad4\beta\leq\alpha.
\end{equation}
Since $\bigcup_{n\geq1}B_{E_n}$ is dense in $B_X$
by \eqref{eq:ball-density},
$\|f\|_{(X_{\R})^*}=1$, and the spaces $E_n$ are increasing, it follows that 
we may choose an integer $n_0\geq1$ sufficiently large and
vectors $z_{n_0}\in S_{E_{n_0}}$ and $y_0\in B_{E_{n_0}}$
such that
\begin{equation}\label{eq:stage-choice}
 2^{-n_0}<\beta,\qquad
 \|z-z_{n_0}\|_X<\varepsilon,\qquad
 f(y_0)>1-\beta.
\end{equation}
Choose $y\in A_{n_0}\subseteq B_{E_{n_0}}$ with
\begin{equation}\label{eq:net-point}
 \|y-y_0\|_{E_{n_0}}
 \overset{\eqref{eq:net-cover}}{\leq}2^{-n_0}
 \overset{    \eqref{eq:stage-choice}}{<} \beta .
\end{equation}
Then
\begin{align}
 f(y)&\quad  =f(y_0)+f(y-y_0)
 \geq f(y_0)-\|f\|_{(X_{\R})^*}\|y-y_0\|_X
 \notag\\
 &  \overset{\eqref{eq:stage-choice},
     \eqref{eq:stage-isometry}}{>}
 1-\beta-\|y-y_0\|_{E_{n_0}}
 \overset{\eqref{eq:net-point}}{>}1-2\beta,
 \label{eq:slice-net-bound}
\end{align}
and hence
\begin{align}
 (1-2^{-n_0})f(y)
 \overset{\eqref{eq:stage-choice},
      \eqref{eq:slice-net-bound}}{>}
 (1-\beta)(1-2\beta)
 =1-3\beta+2\beta^2
 >1-4\beta
 \overset{\eqref{eq:beta}}{\geq} 1-\alpha.
 \label{eq:average-in-slice}
\end{align}
For the chosen $n_0$ and $y\in A_{n_0}$, let
$u_{n_0,y,1},\ldots,u_{n_0,y,N_{n_0,y}}\in S_{E_{n_0+1}}$
be the unit vectors obtained in the construction of $E_{n_0+1}$ (see Proposition \ref{extension}).
They satisfy
\begin{align}\label{2.1}
\frac1{N_{n_0,y}}\sum_{k=1}^{N_{n_0,y}}u_{n_0,y,k}
 \overset{\eqref{orthogonality}}{=}
 (1-2^{-n_0})y
\end{align}
and
\begin{align}\label{2.2}
 \max_{\tau\in\{-1,1\}}
 \|z_{n_0}+\tau u_{n_0,y,k}\|_{E_{n_0+1}} \overset{\eqref{orthogonality}}{=} 
 2, 
 \qquad \forall 1\leq k\leq N_{n_0,y}.
\end{align}
Choose $j\in\{1,\ldots,N_{n_0,y}\}$ such that
$ 
 f(u_{n_0,y,j})
=\mathop{\max}\limits_{1\leq k\leq N_{n_0,y}}f(u_{n_0,y,k})
$,
and set $u:=u_{n_0,y,j}\in   S_X$.
Then
\begin{align*}
 f(u)
 & \,\,  =\max_{1\leq k\leq N_{n_0,y}}f(u_{n_0,y,k})
 \geq\frac1{N_{n_0,y}}
       \sum_{k=1}^{N_{n_0,y}}f(u_{n_0,y,k})\\
 &\overset{\eqref{2.1}}{=}
 (1-2^{-n_0})f(y)
 \overset{\eqref{eq:average-in-slice}}{>}1-\alpha.
\end{align*}
Moreover,
\begin{align*}
 \max_{\tau\in\{-1,1\}}\|z+\tau u\|_X
 &\geq\max_{\tau\in\{-1,1\}}\|z_{n_0}+\tau u\|_X
       -\|z-z_{n_0}\|_X\\
 &\overset{ \eqref{eq:stage-isometry}}{=}
\max_{\tau\in\{-1,1\}} \|z_{n_0}+\tau u\|_{E_{n_0+1}}  -\|z-z_{n_0}\|_X
\\ & \overset{
            \eqref{2.2}}{=} 
 2-\|z-z_{n_0}\|_X
 \overset{\eqref{eq:stage-choice}}{>}2-\varepsilon.
\end{align*}
By Proposition \ref{Prop2.1}, 
$X_{\R}$ has the ADP.

\medskip\noindent\textbf{Step 3: $X$ does not have the DP.}

Define
\begin{align}\label{eq:defT}
T: X &\notag \to X\\ 
v& \mapsto 
 T(v)=\phi(v)x.
\end{align}
Then $T(x)\overset{\eqref{eq:limit-phi}}{=}-x$ and
\begin{equation}\label{eq:T-norm}
 1=\|x\|_X=\|T(x)\|_X\leq\|T\|
 =\|\phi\|_{X^*}\|x\|_X
 \overset{\eqref{eq:limit-phi}}{=}1.
\end{equation}
Fix $u\in B_X$. Note that $|\phi(u)|\in [0,1]$.
If $|\phi(u)|=0$, then
\[
 \|(\Id+T)u\|_X
 \overset{\eqref{eq:defT}}{=}\|u\|_X\leq1.
\]
If $|\phi(u)|>0$, set
\begin{equation}\label{eq:rotation}
 \lambda=\frac{\overline{\phi(u)}}{|\phi(u)|}\in\T,
 \quad 
 \text{then } \quad \lambda u\in B_X,\qquad\phi(\lambda u)=|\phi(u)|.
\end{equation}
Then
\begin{align*}
 \|(\Id+T)u\|_X
 &\quad =\|\lambda(\Id+T)u\|_X
 \overset{\eqref{eq:defT}, \eqref{eq:rotation}}{=}
 \|\lambda u + |\phi(u)| x\|_X
 \\
 & \quad =\||\phi(u)| (\lambda u +x)+(1-|\phi(u)|) \lambda u \|_X\\
 & \quad \leq |\phi(u)| \|\lambda u +x\|_X+(1-|\phi(u)| )\| \lambda u \|_X\\
 &\overset{\eqref{eq:limit-invariant}, \eqref{eq:rotation}}{\leq}
 |\phi(u)| \left(\frac{13}{5}- |\phi(u)| \right)+(1-|\phi(u)|)\\
 & \quad =1+\frac85  |\phi(u)|- |\phi(u)|^2
 =\frac{41}{25}-\left(|\phi(u)|-\frac45\right)^2
 \leq\frac{41}{25}.
\end{align*}
Consequently,
\begin{equation*}
 \|\Id+T\|
 =\sup_{u\in B_X}\|(\Id+T)u\|_X
 \leq\frac{41}{25}<2
 \overset{\eqref{eq:T-norm}}{=}1+\|T\|.
\end{equation*}
\end{proof}

\section*{AI Declaration} 
The authors used ChatGPT 6.1 Sol   to assist with the
discussion and revision of mathematical arguments and
to improve the  language  of the manuscript.  All incorporated
suggestions were critically reviewed, and the mathematical
content was independently verified by the authors, who
take full responsibility for the final manuscript.

\end{document}